\documentclass[11pt,leqno]{article}
\usepackage[margin=1in]{geometry} 
\usepackage{amssymb,amsfonts,amsmath,bbm,mathrsfs,stmaryrd,mathtools}
\usepackage{xcolor}
\usepackage{url}

\usepackage{graphicx}

\usepackage{accents}

\usepackage{extarrows}

\usepackage[shortlabels]{enumitem}
\usepackage{tensor}

\usepackage{xr}
\usepackage[T1]{fontenc}
\usepackage[utf8]{inputenc}

\usepackage[colorlinks,
linkcolor=black!75!red,
citecolor=blue,
pdftitle={},
pdfproducer={pdfLaTeX},
pdfpagemode=None,
bookmarksopen=true,
bookmarksnumbered=true,
backref=page]{hyperref}

\usepackage{tikz}
\usetikzlibrary{arrows,calc,decorations.pathreplacing,decorations.markings,decorations.shapes,intersections,shapes.geometric,through,fit,shapes.symbols,positioning,decorations.pathmorphing}

\makeatletter
\newlength\zig@L
\newlength\zig@La
\newlength\zig@Lb

\newcommand{\xzigrightarrow}[2][]{%
  \mathrel{%
    \settowidth{\zig@La}{$\scriptstyle #2$}%
    \settowidth{\zig@Lb}{$\scriptstyle #1$}%
    \zig@L=\zig@La\relax
    \ifdim\zig@Lb>\zig@L \zig@L=\zig@Lb\fi
    \advance\zig@L by 2.2em\relax
    \tikz[baseline=-0.65ex]{%
      \draw[->,
            line cap=round,
            decorate,
            decoration={zigzag,segment length=4pt,amplitude=1.1pt}]%
        (0,0) -- (\zig@L,0)
        node[midway,above=2pt] {$\scriptstyle #2$}%
        \if\relax\detokenize{#1}\relax\else
          node[midway,below=2pt] {$\scriptstyle #1$}%
        \fi
      ;
    }%
  }%
}
\makeatother

\makeatletter
\newcommand{\squigjoin}{1mu} 

\def\sqleft@{\sim}                    
\def\sqmid@{\sim\mkern-\squigjoin}    

\def\rightsquigarrowfill@{%
  \arrowfill@{\sqleft@}{\sqmid@}{\mkern-4mu\succ}%
}

\newcommand{\xrightsquigarrow}[2][]{%
  \ext@arrow 0359\rightsquigarrowfill@{#1}{#2}%
}
\makeatother

\makeatletter
\newcommand*\circled[1]{\tikz[baseline=(char.base)]{
    \node[shape=circle, draw, inner sep=0pt, 
    minimum height={\f@size},] (char) {\vphantom{WAH1g}#1};}}
\makeatother

\makeatletter
\DeclareRobustCommand\widecheck[1]{{\mathpalette\@widecheck{#1}}}
\def\@widecheck#1#2{%
    \setbox\z@\hbox{\m@th$#1#2$}%
    \setbox\tw@\hbox{\m@th$#1%
       \widehat{%
          \vrule\@width\z@\@height\ht\z@
          \vrule\@height\z@\@width\wd\z@}$}%
    \dp\tw@-\ht\z@
    \@tempdima\ht\z@ \advance\@tempdima2\ht\tw@ \divide\@tempdima\thr@@
    \setbox\tw@\hbox{%
       \raise\@tempdima\hbox{\scalebox{1}[-1]{\lower\@tempdima\box
\tw@}}}%
    {\ooalign{\box\tw@ \cr \box\z@}}}
\makeatother

\usepackage{braket}

\usepackage[amsmath,thmmarks,hyperref]{ntheorem}
\usepackage{cleveref}

\newcommand\nthalias[1]{\AddToHook{env/#1/begin}{\crefalias{lemma}{#1}}}

\nthalias{definition}
\nthalias{example}
\nthalias{examples}
\nthalias{remark}
\nthalias{remarks}
\nthalias{convention}
\nthalias{notation}
\nthalias{construction}
\nthalias{sketch}
\nthalias{theoremN}
\nthalias{propositionN}
\nthalias{corollaryN}
\nthalias{lemma}
\nthalias{proposition}
\nthalias{corollary}
\nthalias{theorem}
\nthalias{conjecture}
\nthalias{question}
\nthalias{assumption}

\creflabelformat{enumi}{#2#1#3}

\crefname{section}{Section}{Sections}
\crefformat{section}{#2Section~#1#3} 
\Crefformat{section}{#2Section~#1#3} 

\crefname{subsection}{\S}{\S\S}
\AtBeginDocument{%
  \crefformat{subsection}{#2\S#1#3}%
  \Crefformat{subsection}{#2\S#1#3}%
}

\crefname{subsubsection}{\S}{\S\S}
\AtBeginDocument{%
  \crefformat{subsubsection}{#2\S#1#3}%
  \Crefformat{subsubsection}{#2\S#1#3}%
}

\theoremstyle{plain}

\newtheorem{lemma}{Lemma}[section]
\newtheorem{proposition}[lemma]{Proposition}

\newtheorem{theorem}[lemma]{Theorem}

\theoremstyle{plain}
\theoremnumbering{Alph}

\theoremstyle{plain}
\theorembodyfont{\upshape}
\theoremsymbol{\ensuremath{\blacklozenge}}

\newtheorem{definition}[lemma]{Definition}
\newtheorem{example}[lemma]{Example}

\newtheorem{remarks}[lemma]{Remarks}

\crefname{definition}{definition}{definitions}
\crefformat{definition}{#2definition~#1#3} 
\Crefformat{definition}{#2Definition~#1#3} 

\crefname{ex}{example}{examples}
\crefformat{example}{#2example~#1#3} 
\Crefformat{example}{#2Example~#1#3} 

\crefname{exs}{example}{examples}
\crefformat{examples}{#2example~#1#3} 
\Crefformat{examples}{#2Example~#1#3} 

\crefname{remark}{remark}{remarks}
\crefformat{remark}{#2remark~#1#3} 
\Crefformat{remark}{#2Remark~#1#3} 

\crefname{remarks}{remark}{remarks}
\crefformat{remarks}{#2remark~#1#3} 
\Crefformat{remarks}{#2Remark~#1#3} 

\crefname{convention}{convention}{conventions}
\crefformat{convention}{#2convention~#1#3} 
\Crefformat{convention}{#2Convention~#1#3} 

\crefname{notation}{notation}{notations}
\crefformat{notation}{#2notation~#1#3} 
\Crefformat{notation}{#2Notation~#1#3} 

\crefname{table}{table}{tables}
\crefformat{table}{#2table~#1#3} 
\Crefformat{table}{#2Table~#1#3}

\crefname{lemma}{lemma}{lemmas}
\crefformat{lemma}{#2lemma~#1#3} 
\Crefformat{lemma}{#2Lemma~#1#3} 

\crefname{proposition}{proposition}{propositions}
\crefformat{proposition}{#2proposition~#1#3} 
\Crefformat{proposition}{#2Proposition~#1#3} 

\crefname{propositionN}{proposition}{propositions}
\crefformat{propositionN}{#2proposition~#1#3} 
\Crefformat{propositionN}{#2Proposition~#1#3} 

\crefname{corollary}{corollary}{corollaries}
\crefformat{corollary}{#2corollary~#1#3} 
\Crefformat{corollary}{#2Corollary~#1#3} 

\crefname{corollaryN}{corollary}{corollaries}
\crefformat{corollaryN}{#2corollary~#1#3} 
\Crefformat{corollaryN}{#2Corollary~#1#3} 

\crefname{theorem}{theorem}{theorems}
\crefformat{theorem}{#2theorem~#1#3} 
\Crefformat{theorem}{#2Theorem~#1#3} 

\crefname{theoremN}{theorem}{theorems}
\crefformat{theoremN}{#2theorem~#1#3} 
\Crefformat{theoremN}{#2Theorem~#1#3} 

\crefname{enumi}{}{}
\crefformat{enumi}{#2#1#3}
\Crefformat{enumi}{#2#1#3}

\crefname{assumption}{assumption}{Assumptions}
\crefformat{assumption}{#2assumption~#1#3} 
\Crefformat{assumption}{#2Assumption~#1#3} 

\crefname{construction}{construction}{Constructions}
\crefformat{construction}{#2construction~#1#3} 
\Crefformat{construction}{#2Construction~#1#3} 

\crefname{sketch}{sketch}{Sketches}
\crefformat{sketch}{#2sketch~#1#3} 
\Crefformat{sketch}{#2Sketch~#1#3} 

\crefname{question}{question}{Questions}
\crefformat{question}{#2question~#1#3} 
\Crefformat{question}{#2Question~#1#3} 

\crefname{equation}{}{}
\crefformat{equation}{(#2#1#3)} 
\Crefformat{equation}{(#2#1#3)}

\numberwithin{equation}{section}

\theoremstyle{nonumberplain}
\theoremsymbol{\ensuremath{\blacksquare}}

\newtheorem{proof}{Proof}
\newcommand\pf[1]{\newtheorem{#1}{Proof of \Cref{#1}}}

\newcommand\bZ{{\mathbb Z}}

\newcommand\cB{{\mathcal B}}

\newcommand\cE{{\mathcal E}}

\newcommand\cL{{\mathcal L}}

\newcommand\cP{{\mathcal P}}

\newcommand\cS{{\mathcal S}}
\newcommand\cT{{\mathcal T}}
\newcommand\cU{{\mathcal U}}

\newcommand\cW{{\mathcal W}}

\DeclareMathOperator{\Aut}{\mathrm{Aut}}

\newcommand{\comment}[1]{}

\title{Power topologies and uniformities on bistochastic operators}
\author{Alexandru Chirvasitu}

\begin{document}

\date{}

\newcommand{\Addresses}{{
  \bigskip
  \footnotesize

  \textsc{Department of Mathematics, University at Buffalo}
  \par\nopagebreak
  \textsc{Buffalo, NY 14260-2900, USA}  
  \par\nopagebreak
  \textit{E-mail address}: \texttt{achirvas@buffalo.edu}

}}

\maketitle

\begin{abstract}
  We prove that the set of invertible Koopman operators attached to the usual interval probability space is dense (and hence dense $G_{\delta}$) in the space of bistochastic operators in the polynomial weak topology. This extends the non-polynomial analogue (due to Brown / Vershik / Ryzhikov / Eisner). We also show that the plain and polynomial weak topologies are distinct on bistochastic operators, confirming the distinction between the two results. Furthermore, the power topologies on the space of bistochastic operators recording weak convergence in powers up to $k+1$, $k\in \mathbb{Z}_{\ge 0}$ are all distinct, as are the analogous power uniformities on the automorphism group of the standard non-atomic Lebesgue probability measure. 
\end{abstract}


\vspace{.5cm}

\noindent{MSC 2020: 47B02; 37A30; 54E15; 54E52; 28D05; 47H09

  
}


\section*{Introduction}

Consider the Lebesgue-measured unit interval $(I:=[0,1],\mu)$ and the \emph{Koopman representation} \cite[Definition 3.8]{gls_erg-join_2003}
\begin{equation*}
  \Aut(\mu)\ni
  T
  \xmapsto{\quad}
  (\bullet)\circ T^{-1}
  \in \cL(L^2(\mu))
  ,\quad
  \cL:=\text{bounded operators}
\end{equation*}
of the group of invertible measure-preserving transformations. The present note is motivated in part by the description of its image's closure in the weak operator topology on the unit ball $\cL(L^2(\mu))_{\le 1}$ as precisely the space $\cB\cS\left(L^2(\mu)\right)$ of \emph{bistochastic} (sometimes \emph{Markov}) operators: positivity-preserving operators which fix the constant function $1\in L^2$ along with their adjoints (e.g. \cite[Theorem 1]{MR192552}, \cite[Remark post Lemma 5.5]{zbMATH01094710} or \cite[Lemma 2.1]{MR5146022}).

Coupled with general interest in the operator-theory literature \cite{MR3101250,zbMATH03749579,zbMATH03132193} in the \emph{polynomial} weak/$\sigma$-weak topology, whereby nets $(T_{\lambda})_{\lambda}$ converge when arbitrary power nets $\left(T^n_{\lambda}\right)_{\lambda}$, $n\in \bZ_{\ge 0}$ do ($\sigma$-)weakly, this raises the natural question of determining the polynomial weak closure of the Koopman image. The main result is that that density persists: 

\begin{theorem}\label{th:bistoc.approx}
  The polynomial weak closure of the set of invertible Koopman operators on $\cL\left(L^2(\mu)\right)$ is precisely the space $\cB\cS\left(L^2(\mu)\right)$ of bistochastic operators on $L^2(\mu)$.
\end{theorem}

\emph{Invertible} Koopman operators are here to be distinguished from the more general operators acting on indicator functions as $\chi_E\mapsto \chi_{T^{-1}E}$ for measure-preserving possibly non-invertible $T$.

It becomes natural at this stage to ask whether the power-topology density is not vacuously a given: distinguishing between the weak and polynomial weak topologies on $\cB\cS\left(L^2(\mu)\right)$ is not a given. That they do indeed differ constitutes a second branch of the note.

\begin{theorem}\label{th:w.not.pw.on.bs}
  The weak and polynomial weak topologies on $\cB\cS\left(L^2\right)$ disagree. 
\end{theorem}

The latter result nods at a broader power-separation phenomenon in evidence when comparing either uniform structures on $\Aut(\mu)$ or companion topologies on $\cB\cS$. Specifically:

\begin{theorem}\label{th:diff.segm}
  Let $k\in \bZ_{\ge 0}$, and consider the uniform structure on $\Aut(\mu)$ generated by the entourages
  \begin{equation}\label{eq:l.unif}
    \left\{(S,T)\in \Aut(I,\mu)^2\ :\ \mu(SA\Delta TA)<\varepsilon\right\}
    ,\quad
    \text{Borel }A\subseteq I
    ,\quad
    \varepsilon>0
  \end{equation}
  with ``$\Delta$'' denoting symmetric set-theoretic difference.

  \begin{enumerate}[(1),wide]
  \item The uniform structures on $\Aut(\mu)$ pulled back from \Cref{eq:l.unif} along
    \begin{equation*}
      \Aut(\mu)
      \ni T
      \xrightarrow{\quad}
      \left(T,\ T^2,\ \cdots,\ T^{k+1}\right)
      \in
      \Aut(\mu)^{k+1}
    \end{equation*}
    are all distinct. 

  \item The same goes for the analogous power topologies on $\cB\cS\left(L^2\right)$.
  \end{enumerate}
\end{theorem}



\section{Polynomial weak approximation of bistochastic operators}\label{se:bistoc}

\pf{th:bistoc.approx}
\begin{th:bistoc.approx}
  The latter space is weakly closed and hence also polynomially so, leaving density as the outstanding issue. For that purpose, the usual partition-based approach to the analogous statement concerning the plain weak topology can be spun into the stronger version.
  
  Functions constant along the atoms of a partition of $I$ into equal-size intervals norm-approximate $L^2$ functions arbitrarily well, so the operators of the form $PMP$ with $M$ bistochastic and
  \begin{equation}\label{eq:part.i}
    P:=\text{orthogonal projection on }\mathrm{span}\left\{\chi_{I_j}\right\}
    ,\quad
    I=\bigsqcup_{j=1}^n I_j
    ,\quad
    \ell\left(I_j\right)=\frac 1n
  \end{equation}
  for varying $n$ constitute a strongly dense subset of $\{M\}$; multiplication being jointly strongly continuous (as in \cite[Solution 113(a)]{hal_hspb_2e_1982}, say), it is enough to approximate $PMP$ for some fixed but arbitrary $n$. 
  
  Partitioning $I$ as in \Cref{eq:part.i}, suppose the sought-after approximation is
  \begin{equation*}
    U^k\ \overset{\text{weak}}{\sim}\ M^k
    ,\quad
    \begin{gathered}
      0\le k\le m\\
      U:=\left(-\circ T\right)\text{ invertible Koopman}\\
      M\text{ bistochastic}.
    \end{gathered}
  \end{equation*}
  The rescaled $\Braket{\chi_{I_i}\mid M\chi_{I_j}}$ constitute a bistochastic matrix $\left(a_{ij}\right)_{1\le i,j\le n}$. We will construct the desired invertible transformation $T\circlearrowright (I,\mu)$ by first partitioning the $I_j$ and then describing the action on the atoms of the partition.

  For arbitrary
  \begin{equation*}
    \mathbf{i}:=
    \left(1\le i_0,\ \cdots,\ i_{m}\le n\right)
    \in
    \left\{1..n\right\}^{m+1}
  \end{equation*}
  let $I_{\mathbf{i}}\subseteq I$ be a subinterval of length $a_{i_0, i_1}a_{i_1, i_2}\cdots a_{i_{m-1},i_m}\ell(I_{j})$; the $I_{\mathbf i}$ partition $I$ upon placing them along it lexicographically. Note that upon extending the notation $I_{\mathbf{i}}$ to possibly shorter tuples $\mathbf{i}\in \left\{1..n\right\}^{m'+1}$ for $0\le m'\le m$, we have sub-partitions 
  \begin{equation*}
    I_{\mathbf{i}'}
    =
    \bigsqcup_{\mathbf{i}\text{ extends }\mathbf{i'}\text{ rightward}}
    I_{\mathbf{i}}
    \quad
    \left(
      \text{using the bistochastic property}
    \right).
  \end{equation*}
  The intention is to have $T^k$, $0\le k\le m$ map $I_{\mathbf{i}}$, $\mathbf{i}\in \left\{1..n\right\}^{m+1}$ to $I_{i_k}$. To that end, simply have $T$ translate every $I_{i_0\cdots i_m}$ into $I_{i_1\cdots i_m}$ so that
  \begin{equation*}
    I_{i_1\cdots i_m}
    =
    \bigsqcup_{1\le j\le n}
    T I_{j\ i_1\cdots i_m}
  \end{equation*}
  lexicographically. This is consistent, again by the bistochastic property. 
\end{th:bistoc.approx}

\begin{remarks}\label{res:same.unif}
  \begin{enumerate}[(1),wide]
  \item\label{item:res:same.unif:unif.det.dns} Recall (e.g. \cite[Problem 20]{hal_hspb_2e_1982}) that the weak and strong topologies of $\cL(H)$ agree on isometries. The corresponding \emph{uniformities} \cite[Definition 7.1]{james_unif_1999} do not (\Cref{ex:unif.punif}); were there such an agreement, the coincidence of the weak and polynomial weak topologies on $\cB\cS:=\cB\cS\left(L^2(\mu)\right)$ would follow from \Cref{th:bistoc.approx} by the following general observation:
    \begin{itemize}[wide]
    \item given uniformities $\cU_0\preceq \cU_1$ (i.e. the latter is \emph{finer}) on a space $X$;
    \item and a $\cU_1$-dense $Y\subseteq X$;
    \item the equality $\cU_0|_Y=\cU_1|_Y$ entails $\cU_0=\cU_1$.
    \end{itemize}
    Or: a uniformity is uniquely determined by its shadow on a dense subspace; this follows, for instance, from \cite[Proposition 8.15]{james_unif_1999}.

  \item\label{item:res:same.unif:diff.top.bs} In the present context there is a converse to the preceding observation: the agreement of the weak$^*$ and polynomial weak$^*$ topologies on $\cB\cS\left(L^2:=L^2(\mu)\right)$ is \emph{equivalent} to the agreement of the two uniform structures on the space of invertible Koopman operators. Indeed, on examining invertible-Koopman nets
    \begin{equation*}
      \left(T_{\lambda}\right)_{\lambda}
      \quad\text{and}\quad
      \left(S_{\lambda}\right)_{\lambda}
      ,\quad
      T_{\lambda}\overset{\quad\text{uniformly close}\quad}{\quad\sim\quad}S_{\lambda},
    \end{equation*}    
    one can always assume $\left(T_{\lambda}\right)$ convergent in $\cB\cS\left(L^2\right)$ by the weak$^*$ compactness of the latter space.

  \item\label{item:res:same.unif:gen.diff.top.bs} The preceding item illustrates a more general principle useful in promoting uniformity-to-topology distinction: given nets mutually approaching in a weak uniformity $\cU_0$ but separated by a stronger $\cU_1$, the induced topologies $\cT_{0,1}$ on a space containing $X$ $\cU_1$-densely are also distinct provided the nets are $\cT_0$-convergent.
  \end{enumerate}  
\end{remarks}

The weak and polynomial weak uniformities do not coincide on the unitary group $U(\ell^2)$: 

\begin{example}\label{ex:unif.punif}
  There are unitary operators $U_{0,n}$ and $U_{1,n}$, $n\in \bZ_{\ge 0}$ with
  \begin{equation*}
    U_{0,n},\ U^2_{0,n},\ U_{1,n},\
    \xrightarrow[\quad n\quad]{\quad\text{weakly}\quad}
    0
    \quad\text{and}\quad
    U^2_{1,n}
    \xrightarrow[\quad n\quad]{\quad\text{weakly}\quad}
    1;
  \end{equation*}
  this renders $U_{i,n}$ uniformly close for $n\gg 0$, while $U^2_{i,n}$ are not so. As to the existence claim, the $n^{th}$ power of the \emph{bilateral shift} \cite[Problem 84]{hal_hspb_2e_1982} will do for  $U_{0,n}$, while $U_{1,n}$ can be chosen so as to interchange standard basis elements $e_{k}$ and $e_{k'}$ for some partition of the basis indices into pairs $\{k,k'\}$ with all $|k-k'|\ge n$.
\end{example}

The kind of uniformity disagreement visible in \Cref{ex:unif.punif} also obtains at the measure-algebra-automorphism level: it is possible for $S,T\in \Aut(I,\mu)$ to be arbitrarily close while their squares are uniformly far with respect to (the \emph{left} uniformity \cite[\S 7.1, post Definition 7.3]{james_unif_1999}, generated by) \Cref{eq:l.unif}. 

\begin{remarks}\label{res:cmpr.unifs}
  \begin{enumerate}[(1),wide]
  \item The left uniformity generated by \Cref{eq:l.unif} cannot agree with that inherited from the weak operator topology: the latter, in turn generated by
    \begin{equation}\label{eq:wot.unif}
      \left\{(S,T)\in \Aut(I,\mu)^2\ :\ \left|\mu(A\cap SB)-\mu(A\cap TB)\right|<\varepsilon\right\}
      ,\quad
      \text{Borel }
      A,B\subseteq I
      ,\quad
      \varepsilon>0,
    \end{equation}
    is inversion-invariant, whereas the former is not. Indeed, \cite[pp. 63-64]{halm_erg_1960} argues that the left/right uniformities are \emph{incomplete} in the sense of \cite[Definition 12.1]{james_unif_1999}, while the bilateral uniformity generated by
    \begin{equation}\label{eq:bil.unif}
      \left\{(S,T)\in \Aut(I,\mu)^2
        \ :\
        \begin{gathered}
          \mu(SA\Delta TA)<\varepsilon\\
          \mu(S^{-1}A\Delta T^{-1}A)<\varepsilon
        \end{gathered}
      \right\}      
    ,\quad
    \text{Borel }A\subseteq I
    ,\quad
    \varepsilon>0
  \end{equation}
  is.

  \item Nor can the (inversion-invariant) \Cref{eq:wot.unif} and \Cref{eq:bil.unif} agree: the latter is, as just noted, complete, whereas the former cannot be by \Cref{th:bistoc.approx}. See also \cite[Theorem 1.1]{MR2957245} for further elaboration on \Cref{eq:wot.unif}. 
  \end{enumerate}
\end{remarks}

It is a natural question, at this point, whether the weak and polynomial weak topologies agree on $\cB\cS$: were this the case, \Cref{th:bistoc.approx} would be a consequence of \cite[Lemma 2.1]{MR5146022}. The following result confirms that this is not the case.

\begin{proposition}\label{pr:w.w2.unif.dist}
  The uniformities on invertible Koopman operators induced by the weak operator topology and, respectively, its restriction along $T\xmapsto{}(T,T^2)$ are distinct. 
\end{proposition}
\begin{proof}
    It will suffice, for arbitrary $m\in \bZ_{> 0}$, to produce invertible measure-preserving $S,T$ on $I:=[0,1]$ (equipped with its Lebesgue measure $\mu$, as before) so that
  \begin{enumerate}[(a),wide]
  \item\label{item:pr:w.w2.unif.dist:agree.2m} we have agreement 
    \begin{equation*}
      \forall\left(0\le k< 2^{m}\right)
      \left(S I_{m,k} = T I_{m,k}\right)
      ,\quad
      I_{m,k}:=\left[\frac{k}{2^m},\frac{k+1}{2^m}\right];
    \end{equation*}

  \item\label{item:pr:w.w2.unif.dist:st.l12} with
    \begin{equation*}
      \forall\left(0\le k<2^m\right)
      \left(
        I_{m+1,k}
        \xrightarrow[\quad\substack{\text{increasingly}\\\text{linearly}}\quad]{\quad S,T\quad}
        I_{m+1,2k}
      \right);
    \end{equation*}

  \item\label{item:pr:w.w2.unif.dist:s} while
    \begin{equation*}
      \forall\left(0\le j<2^{m-1}\right)
      \left(
        \begin{aligned}
          I_{m+1,2^m+2j}
          &\xrightarrow[\quad\substack{\text{increasingly}\\\text{linearly}}\quad]{\quad S\quad}
          I_{m+1,2j+1}\\
          I_{m+1,2^m+2j+1}
          &\xrightarrow[\quad\substack{\text{increasingly}\\\text{linearly}}\quad]{\quad S\quad}
          I_{m+1,2(2^{m}-1-j)+1}
        \end{aligned}        
      \right)
    \end{equation*}

  \item\label{item:pr:w.w2.unif.dist:t} and    
    \begin{equation*}
      \forall\left(0\le j<2^{m-1}\right)
      \left(
        T I_{m+1,2^m+2j}
        =
        S I_{m+1,2^m+2j+1}
        \ \wedge\ 
        T I_{m+1,2^m+2j+1}
        =
        S I_{m+1,2^m+2j}
      \right),
    \end{equation*}
    universal quantifiers ranging over $k$ and $j$ for which the symbols involved are defined. 
  \end{enumerate}
  \Cref{item:pr:w.w2.unif.dist:agree.2m} makes $S$ and $T$ arbitrarily \Cref{eq:l.unif}-close for $m\gg 0$, while $S^2$ and $T^2$ are uniformly \Cref{eq:wot.unif}-far no matter how large $m$ is:
  \begin{equation*}
    m\gg 0
    \quad
    \xRightarrow{\quad}
    \quad
    \mu\left(\left[0,\frac 12\right]\cap S^2\left[0,\frac 12\right]\right)=\frac 38
    \quad\text{and}\quad
    \mu\left(\left[0,\frac 12\right]\cap T^2\left[0,\frac 12\right]\right)=\frac 18.
  \end{equation*}
  The conclusion follows from the fact that the former uniformity dominates the latter.
\end{proof}

\pf{th:w.not.pw.on.bs}
\begin{th:w.not.pw.on.bs}
    This follows from \Cref{pr:w.w2.unif.dist} and \Cref{res:same.unif}\Cref{item:res:same.unif:diff.top.bs}. 
\end{th:w.not.pw.on.bs}


\Cref{pr:w.w2.unif.dist} and \Cref{th:w.not.pw.on.bs} suggest a more extensive analysis of the extent to which power maps on $\Aut(I,\mu)$ and $\cB\cS$ are topologically and/or uniformly ``independent''. With that in view, consider the following posets of topologies and uniformities.

\begin{definition}\label{def:pow.tops.unifs}
  For a topology $\cT$ on a semigroup $S$ making the multiplication separately continuous and $\cE\subseteq \bZ_{>0}$, write $\cT_{\cE}$ for the pullback topology induced by 
  \begin{equation*}
    S
    \ni
    s
    \xmapsto{\quad}
    \left(s^n\right)_{n\in \cE}
    \in
    \left(S^{\cE},\ \text{product topology }\cT^{\cE}\right).
  \end{equation*}
  The notation extends by analogy to
  \begin{itemize}[wide]
  \item topologies on \emph{groups} and subsets $\cE\subseteq \bZ_{\ne 0}$ instead;
  \item or uniform structures $\cU$, whereby the corresponding power uniformity would be $\cU_{\cE}$. 
  \end{itemize}
  This applies in particular to
  \begin{itemize}[wide]
  \item the weak topology $\left(\cB\cS\left(L^2\right),\cW\right)$;
  \item the induced (left=right) uniformity $\left(\Aut(\mu),\cW=\cW_{\ell}\right)$ on invertible Koopman operators;
  \item or the natural left uniformity $\left(\Aut(\mu),\cU_{\ell}\right)$ (or its right counterpart) of \cite[pp. 63-64]{halm_erg_1960}. 
  \end{itemize}
\end{definition}

\pf{th:diff.segm}
\begin{th:diff.segm}
  \Cref{th:diff.struct.k.k1} will provide the induction step, and \Cref{pr:w.w2.unif.dist} the base case of $k=0$. Since compactness of the weak topology on the unit ball permits the assumption that the $S=S_m$ and $T=T_m$ produced in \Cref{pr:w.w2.unif.dist} converge weakly in $\cB\cS$, \Cref{th:diff.struct.k.k1}\Cref{item:th:diff.struct.k.k1:unif.top} is operative. Distinguishing uniformities on $\Aut(\mu)$ thus also differentiates topologies on $\cB\cS$ via \Cref{res:same.unif}\Cref{item:res:same.unif:gen.diff.top.bs}. 
\end{th:diff.segm}

In \Cref{th:diff.struct.k.k1}, for nets $x:=(x_{\alpha})$, $y:=(y_{\alpha})$ in a uniform space $(X,\cU)$, write $x\sim_{\cU} y$ if they mutually approach arbitrarily with respect to the uniformity.

\begin{theorem}\label{th:diff.struct.k.k1}
  Let $k\in \bZ_{\ge 0}$.
  
  \begin{enumerate}[(1),wide]
  \item\label{item:th:diff.struct.k.k1:unif} If $\left(S_{\bullet},T_{\bullet}\right)$ are pairs of nets in $\Aut(\mu)$ satisfying the condition
    \begin{equation*}
      \cP_k
      \quad:\quad
      \forall\left(1\le d\le k+1\right)
      \left(
        S^d_{\bullet}
        \sim_{\text{\Cref{eq:l.unif}}}
        T^d_{\bullet}
      \right)
      ,\quad
      S^{k+2}_{\bullet}
      \not\sim_{\text{\Cref{eq:wot.unif}}}
      T^{k+2}_{\bullet}
    \end{equation*}
    then pairs of nets $(S'_{\bullet},T'_{\bullet})$ exist satisfying the analogous $\cP_{k'}$ for
    \begin{equation*}
      k':=k+(t-s)\left\lfloor\frac{k+1}{s}\right\rfloor
      \quad\text{and any}\quad
      1\le s\le t\in \bZ_{>0}.
    \end{equation*}
    
  \item\label{item:th:diff.struct.k.k1:unif.top} If furthermore $S_{\bullet}$ and $T_{\bullet}$ converge with respect to the power topology $\cW_{k}$ on $\cB\cS$ pulled back from the product weak topology along
    \begin{equation*}
      \cB\cS\ni S
      \xmapsto{\quad}
      \left(S,\ S^2,\ \cdots,\ S^{k+1}\right)
      \in \cB\cS^{k+1}
    \end{equation*}
    $S'_{\bullet}$ and $T'_{\bullet}$ can be assumed to converge in the analogous power topology $\cW_{k'}$.
  \end{enumerate}
\end{theorem}
\begin{proof}
  Both statements will follow from a direct construction of $S\mapsto S'$, $S\in \Aut(\mu)$, by direct inspection and elementary computational checks we do not undertake.

  Having fixed $1\le s\le t\in \bZ_{>0}$ throughout, consider the probability space
  \begin{equation*}
    \left(J,\nu\right)
    ,\quad
    \begin{aligned}
      J&:=I\times \bZ/t\\
      \nu&:=\text{scaled Lebesgue, uniformly distributed among the $t$ components}.
    \end{aligned}
  \end{equation*}
  Being \emph{standard} non-atomic, $(J,\nu)\cong (I,\mu)$ \cite[Theorem 2.11]{gls_erg-join_2003} so in particular $\Aut(\nu)\cong \Aut(\mu)$. Given $S\in \Aut(\mu)$, simply take for $S'$ the transformation cycling through the $t$ components of $J$, implementing
  \begin{itemize}[wide]
  \item the identity during the first $t-s$ steps, i.e. in transitioning from $\{0\}\times I$ to $\{t-s\}\times I$;
  \item and $S$ for each of the remaining $s$ steps, transitioning from $\{t-s\}\times I$ back to $\{0\}\times I$. 
  \end{itemize}
  The construction preserves both mutual approach and convergence in the space of bistochastic operators. 
\end{proof}


\addcontentsline{toc}{section}{References}

\def\polhk#1{\setbox0=\hbox{#1}{\ooalign{\hidewidth
  \lower1.5ex\hbox{`}\hidewidth\crcr\unhbox0}}}


\Addresses

\end{document}